%


\magnification=\magstep1
\input amstex
\documentstyle{amsppt}
\def\pcf{\text{pcf}}                      %
\def\cof{\text{cof}}                      %
\def\aloe{\aleph_{\omega}}                %
\NoBlackBoxes
\topmatter
\title
Possible pcf algebras
\endtitle
\author
Thomas Jech and Saharon Shelah
\endauthor
\thanks
The first author was supported in part by an NSF grant
DMS-8918299, and by the U.S.--Israel Binational Science
Foundation.
\endgraf
The second author was partially supported by the U.S.--Israel
Binational Science Foundation. Publication No. 476\endthanks
\affil
The Pennsylvania State University \\
The Hebrew University and Rutgers University
\endaffil
\address
Department of Mathematics, The Pennsylvania State University,
University Park, PA 16802, USA
\endgraf
School of Mathematics, The Hebrew University, Jerusalem, Israel, and
Department of Mathematics, Rutgers University, New Brunswick, NJ 08903, USA
\endaddress
\email jech\@math.psu.edu, shelah\@math.huji.ac.il \endemail
\abstract
There exists a family $\{B_{\alpha}\}_{\alpha<\omega_1}$
of sets of countable ordinals such that
\roster
\item $\max B_{\alpha}=\alpha$,
\item if $\alpha\in B_{\beta}$ then $B_{\alpha}\subseteq B_{\beta}$,
\item if $\lambda\leq \alpha$ and $\lambda$ is a limit ordinal then 
$B_{\alpha}\cap\lambda$ is not in the ideal generated by the $B_{\beta}$, 
$\beta<\alpha$, and by the bounded subsets of $\lambda$,
\item  there is a partition $\{A_n\}_{n=0}^{\infty}$ of $\omega_1$ 
such that for every $\alpha$ and every $n,$ $B_{\alpha}\cap A_n$ is finite.
\endroster
\endabstract
\endtopmatter

\document
\baselineskip 20pt
\subhead
1. Introduction
\endsubhead

In \cite3, \cite4, \cite5 and \cite6 the second
author developed the theory of possible cofinalities (pcf), and
proved, among others, that if $\aloe$ is a strong limit
cardinal then $2^{\aloe}<\aleph_{(2^{\aleph_0})^+}$ as well
as $2^{\aloe}<\aleph_{\omega_4}$.  The latter inequality is
established via an analysis of the structure of pcf; in particular, it
is shown that if $\aleph_4\leq
\left|\pcf\{\aleph_n\}_{n=0}^{\infty}\right|$ then a certain structure
exists on $\omega_4$, and then it is proved that such a structure is
impossible.  (Cf. \cite5, \cite1 and \cite2 for details.)
One might hope that by investigating this structure one could possibly
derive a contradiction for $\aleph_3, \,\aleph_2$ or even $\aleph_1.$

A major open problem in the theory of singular cardinals (or in the
pcf theory) is whether it is consistent that $\aloe$ is strong limit
and $2^{\aloe}>\aleph_{\omega_1}$; or whether the set
$\pcf\,\{\aleph_n\}_{n=1}^{\infty}$ can be uncountable.  If we make this
assumption, we obtain a certain structure on $\omega_1$.  The
structure is described in Theorem 2.1.  Unlike in the $\omega_4$-case,
the structure so obtained is not impossible: in Theorem 3.1 we show
that there exists a structure on $\omega_1$ described in
the abstract, and consequently has the properties given by Theorem 2.1.

In Section 2, all facts on Shelah's pcf theory not proved explicitly
can be found in the expository articles \cite1 and \cite2.  In Section
3 we assume rudimentary knowledge of forcing.

\subhead
2.  A consequence of ``$\pcf\;\{\aleph_n\}_{n=0}^{\infty}$ is uncountable''
\endsubhead

\proclaim
{Theorem 2.1} If $\pcf\;\{\aleph_n\}_{n=0}^{\infty}$ is uncountable,
then there exist sets $B_{\alpha}$, $\alpha<\omega_1$, of countable
ordinals with the following properties:
\roster
\item"(a)"  For every $\alpha<\omega_1$, \;\;$\max\;B_{\alpha}=\alpha$.
\item"(b)"  For all $\alpha,\beta<\omega_1$, if $\alpha\in B_{\beta}$ then 
$B_{\alpha}\subseteq B_{\beta}$.
\item"(c)"  For every limit ordinal $\lambda<\omega_1$, 
$B_{\lambda}\cap\lambda$ is unbounded in $\lambda$.
\item"(d)"  There is a closed unbounded set $C$ of countable limit ordinals 
such that for all $\lambda\in C$ and for all $\alpha\geq\lambda$,  the set 
$B_{\alpha}\cap\lambda$ is not in the ideal generated by the sets $B_{\beta}$, 
$\beta<\alpha$, and by bounded subsets of $\lambda$.  (I.e. 
$B_{\alpha}\cap\lambda \nsubseteq\gamma\cup B_{\beta_1}\cup\cdots\cup 
B_{\beta_k}$, for any $\gamma<\lambda$, and any $\beta_1,\dots,\beta_k<\alpha$.)
\item"(e)"  Every unbounded set $X\subseteq\omega_1$ has an initial segment 
$X\cap\gamma$ that is not in the ideal generated by the sets $B_{\alpha}$, 
$\alpha<\omega_1$.
\item"(f)"  Moreover, (e) remains true in every extension $M$ of the ground 
model that preserves cardinals and cofinalities, and has the property that 
every countable set of ordinals in $M$ is covered by a countable set in the 
ground model.
\endroster
\endproclaim

\demo
{Proof} Let $a=\pcf\,\{\aleph_n\}_{n=0}^{\infty}$ and assume that $a$ is
uncountable.  Applying the pcf theory, one obtains (cf. \cite6, Main
Theorem) sets $b_{\lambda}$, $\lambda\in a$, ({\it generators})
together with sequences of functions $f_i^{\lambda}$ $(i<\lambda)$ in
$\prod a$.  As $a$ contains all regular cardinals
$\lambda<\aleph_{\omega_1}$, we let, for each $\alpha<\omega_1$
$$
B_{\alpha}=\{\xi:\aleph_{\xi+1}\in b_{\aleph_{\alpha+1}}\}.
$$

Property (a) is immediate.  Property (b) is the
{\it transitivity} of generators; such generators can be found
(cf. \cite1, Lemma 6.9).

Property (c) is a consequence of the fact that for every countable
limit ordinal $\lambda$, there exists an increasing sequence
$\alpha_n$, $n<\omega$, with limit $\lambda$, and an ultrafilter $D$
on $\omega$ such that $\cof\left(\overset{\infty}\to{\underset{n=0}
\to\prod}\aleph_{\alpha_n+1}/D\right)=\aleph_{\lambda+1}$
(cf. \cite1, Theorem 2.1).

Property (d): Let $\gamma_i$, $i<\omega_1$, be a continuous increasing
sequence of countable ordinals constructed as follows: Given
$\gamma_i$, we first note that
$\aleph_{\omega_1+1}\in\pcf\,[\aleph_{\gamma_i+1},\;\aleph_{\omega_1})$
(by \cite1, Theorem 2.1), and by the Localization Theorem \cite6,
there is a $\gamma_{i+1}<\omega_1$ such that
$\aleph_{\omega_1+1}\in\pcf\,[\aleph_{\gamma_i+1},\;\aleph_{\gamma_{i+1}})$.
Let $C$ be the set of all limit points of the sequence
$\{\gamma_i\}_{i<\omega_1}$.

Now let $\lambda\in C$, $\alpha\geq\lambda$, $\gamma<\lambda$, and
$\beta_1,\dots,\beta_k<\alpha$.  We find $\gamma_i$ such that
$\gamma<\gamma_i<\lambda$.  By \cite{1}, Theorem 2.1, we have
$\aleph_{\alpha+1}\in \pcf\,[\aleph_{\gamma_i+1},\aleph_{\gamma_{i+1}})$
and so there is an ultrafilter $D$ on $[\gamma_i+1,\gamma_i)$ such
that $\cof\left(\prod\aleph_{\xi+1}/D\right)=\aleph_{\alpha+1}$.  By
the definition of generators, we have $B_{\alpha}\in D$ while
$B_{\beta_i}\notin D$ $(i=1,\dots,k)$, and (d) follows.

Property (e):\quad If $X\subseteq\omega_1$ is unbounded, then
$\max\pcf\,\{\aleph_{\alpha+1}:\alpha\in X\}\geq\aleph_{\omega_1}$, and
by the Localization Theorem, there is a countable $\gamma$ such that
$\max\pcf\,\{\aleph_{\alpha+1}:\alpha\in
X\cap\gamma\}\geq\aleph_{\omega_1}$.  Now if $\alpha_1,\dots,\alpha_k$
are countable ordinals, we cannot have $X\cap\gamma\subseteq
B_{\alpha_1}\cup\cdots\cup B_{\alpha_k}$, because
$\max\pcf\,(b_{\aleph_{\alpha_1}+1}\cup\cdots\cup
b_{\aleph_{\alpha_k}+1})=\max\,\{\aleph_{\alpha_i+1}:i=1,\dots,k\}<
\aleph_{\omega_1}$.

Property (f):\quad Let $M$ be an extension of the ground model $V$
that preserves cardinals and cofinalities, and assume further that
every countable set of ordinals in $M$ is covered by a countable set
in $V$.

To show that (e) is true in $M$, it suffices to show that the
generators $b_{\lambda}$ are generators of the $\pcf$ structure in
$M$.  For that, it is enough to verify that the sequences
$f_i^{\lambda}$ $(i<\lambda)$ are increasing cofinal sequences in $\prod
a$ (modulo the appropriate ideals $J_{<\lambda}$).  Since $M$ has the
same cardinals and cofinalities, the claim follows upon the
observation that for every regular $\lambda<\aleph_{\omega_1}$, every
function $f\in\prod b_{\lambda}$ in $M$ is majorized by some function
$g\in\prod b_{\lambda}$ in $V$.  
\enddemo

\subhead
3.\quad Existence of the family $\{B_{\alpha}\}_{\alpha<\omega_1}$
\endsubhead

\proclaim
{Theorem 3.1} There exist a
partition $\{A_n\}_{n=0}^{\infty}$ of $\omega_1$, and a family
$\{B_{\alpha}\}_{\alpha<\omega_1}$ of countable sets of countable
ordinals such that
\roster
\item"(a)"  For every $\alpha<\omega_1$, $\max B_{\alpha}=\alpha.$
\item"(b)"  For all $\alpha,\beta<\omega_1$, if $\alpha\in B_{\beta}$ then 
$B_{\alpha}\subseteq B_{\beta}$.
\item"(c)"  For every limit ordinal $\lambda<\omega_1$ and for all 
$\alpha\geq\lambda$, $B_{\alpha}\cap\lambda\nsubseteq \gamma\cup 
B_{\beta_1}\cup\cdots\cup B_{\beta_k}$ for any $\gamma<\lambda$ and any 
$\beta_1,\dots,\beta_k<\alpha$.
\item"(d)"  For all $\alpha<\omega_1$ and all $n$, $B_{\alpha}\cap A_n$ 
is finite.
\endroster
\endproclaim

\proclaim
{Corollary 3.2} If $M$ is any $\aleph_1$-preserving
extension of $V$, then every unbounded set $X\subseteq \omega_1$ in
$M$ has an initial segment $X\cap\gamma$ that is not in the ideal
generated by the sets $B_{\alpha}$, $\alpha<\omega_1$.
\endproclaim

\demo
{Proof} By (d), any set in the ideal has a finite intersection with
each $A_n$.  If $X\subseteq\omega_1$ is unbounded then some $X\cap
A_n$ is uncountable, and so some $(X\cap\gamma)\cap A_n$ is infinite.
Hence $X\cap\gamma$ is not in the ideal.  
\enddemo

To construct the structure described in Theorem 3.1 we shall first
define a forcing notion and prove that it forces such a structure to exist
in the generic extension. The forcing notion that we use satisfies the
countable chain condition and consists of finite conditions consisting
of countable ordinals and relations between countable ordinals. Using
a general method due to the second author \cite7 we then conclude that
such a structure exists in $V.$

\proclaim{Definition 3.3} \endproclaim

A forcing condition is a quadruple $p=(S_p,\;\pi_p,\;b_p,\;u_p)$ such that

\roster
\item"(i)"\quad $S_p$ is a finite subset of $\omega_1,$
\item"(ii)"\quad  $b_p$ is a function from $S_p\times S_p$ 
into $\{0,1\}$ such that
$$
\gathered
b_p(\alpha,\alpha)=1\hskip 1.5in (\alpha\in S_p) \\
b_p(\alpha,\beta)=0\hskip 1.5in (\alpha,\beta\in S_p,\quad \alpha<\beta) \\
\text{if }b_p(\alpha,\beta)=1\text{ and }b_p(\beta,\gamma)=1\;\;\text{then }
b_p(\alpha,\gamma)=1\qquad\qquad(\alpha,\beta,\gamma\in S_p)
\endgathered
$$
\item"(iii)"\quad  $u_p$ is a natural number,
\item"(iv)"\quad   $\pi_p$ is a function from $S_p$ into 
$\{0,...,u_p-1\}$ such that for all $\alpha$ and $\beta$ in $S_p$, if
$b_p(\alpha,\beta)=1$ and $\beta <\alpha$ then $\pi_p(\beta)\ne\pi_p(\alpha),$

\endroster
[Motivation: $S$ is the support of the condition, $\pi(\alpha)=n$
forces $\alpha\in A_n$, $b(\alpha,\beta)=1$ forces $\beta\in
B_{\alpha}$ and $b(\alpha,\beta)=0$ forces $\beta\notin B_{\alpha}$.]

A condition $r=(S_r,\;\pi_r,\;b_r,\;u_r)$ is stronger than 
$p=(S_p,\;\pi_p,\;b_p,\;u_p)$ if
\roster
\item"(i)"  $S_r\supseteq S_p,$
\item"(ii)"  $b_r$ extends $b_p,$
\item"(iii)"  $\pi_r$ extends $\pi_p$
\item"(iv)"  $u_r \geq u_p,$
\item"(v)"  for all $\alpha\in S_p$ and all $\beta\in S_r-S_p$, if 
$b_r(\alpha,\beta)=1$ then $\pi_r(\beta)\geq u_p$.
\endroster

It is easy to verify that ``stronger than'' is a transitive relation.

\proclaim {Definition 3.4} \endproclaim If $p=(S_p,\;\pi_p,\;b_p,\;u_p)$
is a condition and $\eta$ is a countable ordinal, we let
$$
p\restriction\eta=(S_p\cap\eta,\;\pi_p\restriction\eta,\;b_p\restriction 
(\eta\times\eta),\;u_p).
$$
Clearly, $p\restriction\eta$ is a condition and $p$ is stronger than
$p\restriction\eta.$

\proclaim {Lemma 3.5 (Ammalgamation)}  If $p$ and $q$ are conditions and
$\eta$ a countable ordinal such that $q$ is stronger than $p\restriction\eta$
and $S_q \subseteq \eta$ then there exists a condition $r$ such that $r$
is stronger than both $p$ and $q$ (and such that $S_r=S_p\cup S_q$).
\endproclaim

\demo{Proof}
Note that $u_q\geq u_p.$ We let $S_r=S_p\cup S_q,$ $\pi_r=\pi_p\cup\pi_q$
and $u_r=u_q.$ We define $b_r$ as follows: if  $\alpha$ and $\beta$ are
both in $S_p$ (both in $S_q$) the we let $b_r(\alpha,\beta)=b_p(\alpha,\beta)$
(we let $b_r(\alpha,\beta)=b_q(\alpha,\beta).)$ If $\alpha\geq\eta$ is in
$S_p$ and if $\beta<\eta$ is in $S_q - S_p$ then we let $b_r(\alpha,\beta)=1$ if
and only if there exists a $\gamma<\eta$ in $S_p$ such that $b_p(\alpha,\gamma)
=1$ and $b_q(\gamma,\beta)=1.$ Otherwise we let $b_r(\alpha,\beta)=0.$

To verify that $r$ is a condition, it is easy to see that condition (ii)
from the definition is satisfied. To verify (iv), the only case we need
to worry about is when $b_r(\alpha,\beta)=1$ where $\alpha\geq\eta$ is
in $S_p$ and $\beta<\eta$ is in $S_q - S_p.$ In this case, $\pi_q(\beta)\ge u_p$
(because $q$ is stronger than $p\restriction\eta$ 
and $b_q(\gamma,\beta)=1$ for some $\gamma\in S_p\cap\eta$) 
while $\pi_p(\alpha)<u_p,$ and so $\pi_r(\beta)\ne\pi_r(\alpha).$

Since $r\restriction\eta=q$, $r$ is stronger than $q.$ In order to show that
$r$ is stronger than $p$ we only need to verify condition (v), and only for
the case when $\alpha\geq\eta$ is in $S_p$ and $\beta<\eta$ is in $S_q - S_p.$
This is however exactly the argument in the preceding paragraph.
\enddemo

\proclaim{Lemma 3.6} The forcing satisfies the countable chain condition.
\endproclaim

\demo{Proof} Given $\aleph_1$ conditions, we first find $\aleph_1$ of them
whose supports form a $\Delta$-system, with a root $A$, i.e. $S_{p_\xi}\cap
S_{p_\eta}=A$ whenever $\xi<\eta$, and such that $\beta<\alpha$ whenever
$\beta\in S_{p_\xi}$ and $\alpha\in S_{p_\eta}-A$.  Then $\aleph_1$ of them have the
same restrictions of $\pi$ and $b$ to the root $A,$ and the same $u$.

Now it follows from Lemma 3.5 that any two such conditions are compatible. 
\enddemo

Let $G$ be a generic set of conditions.  In $V[G]$, we let, for each
$\alpha<\omega_1$ and each $n<\omega$,
$$
B_{\alpha}=\{\beta:b(\alpha,\beta)=1 
\text{ for some condition } (S,\;\pi,\;b,\;u)\in G\},
\tag3.7
$$
$$
A_n=\{\alpha:\pi(\alpha)=n \text{ for some condition } (S,\;\pi,\;b,\;u)\in G\}.
\tag3.8
$$

Clearly,  $\max B_{\alpha}=\alpha$, and if
$\alpha\in B_{\beta}$ then $B_{\alpha}\subseteq B_{\beta}$.
The sets $A_n$ are mutually disjoint subsets of $\omega_1.$

\proclaim{Lemma 3.9} For every $\alpha<\omega_1$ the set of all conditions
$p$ with $\alpha\in S_p$ is dense. For every $n$ the set of all conditions
$p$ with $u_p\geq n$ is dense.
\endproclaim

\demo{Proof} If $q$ is a condition and $\alpha\notin S_q$ then let $S_p=
S_q\cup\{\alpha\},$ let $b_p(\alpha,\alpha)=1,$ $u_p=u_q+1$ and $\pi_p(\alpha)=
u_q.$ Then $p$ is a condition stronger than $q.$ The proof of the
second statement is similar.
\enddemo

\proclaim{Corollary 3.10}
$\{A_n\}_{n=0}^{\infty}$ is a partition of $\omega_1$.
\endproclaim

\proclaim{Lemma 3.11} For all $\alpha<\omega_1$ and all $n$,
$B_{\alpha}\cap A_n$ is finite.
\endproclaim

\demo
{Proof} Let $\alpha$ and $n$ be given, and let
$p=(S_p,\;\pi_p,\;b_p,\;u_p)$ be a condition.  We shall find a stronger
condition $q$ that forces that $B_{\alpha}\cap A_n$ is finite.

There is a condition $q=(S_q,\;\pi_q,\;b_q,\;u_q)$ stronger than $p$
such that $\alpha\in S_q$ and that $u_q>n.$ We claim that $q$ forces that
$B_{\alpha}\cap A_n\subseteq S_q$.

If $\beta$ is an ordinal not in $S_q$ and if $r=(S_r,\;\pi_r,\;b_r,\;u_r)$ is
a stronger condition that forces $\beta\in B_{\alpha}$ then because
$b_r(\alpha,\beta)=1$, we have $\pi_r(\beta)\geq u_q>n$, and so $r$
forces $\beta\notin A_n$.  Thus $q$ forces $B_{\alpha}\cap A_n\subseteq S_q$.
\enddemo

\proclaim{Lemma 3.12} Let $\lambda<\omega_1$ be a limit
ordinal, let $\alpha\geq\lambda$, and let $\gamma<\lambda$ and
$\alpha_1,\dots,\alpha_k<\alpha$.  There exists a $\beta\geq\gamma$,
$\beta<\lambda$, such that $\beta\in B_{\alpha}$ and $\beta\notin~B_{\alpha_1},
\dots,\beta\notin B_{\alpha_k}$.
\endproclaim

\demo
{Proof} Let $p=(S_p,\;\pi_p,\;b_p,\;u_p)$ be a condition.  We may assume
that $\alpha,\alpha_1,\dots,\alpha_k\in S_p$.  Let
$\beta<\lambda$ be such that $\beta\geq\gamma$ and $\beta\notin S_p$.

Let $\eta=\alpha+1$ and $S=S_p \cap\eta.$ We let $S_q=S\cup\{\beta\}$, 
$u_q=u_p+1,$ $\pi_q \restriction S = \pi_p \restriction S,$
$\pi_q(\beta)=u_p,$ $b_q \restriction (S \times S) = b_p \restriction
(S \times S),$ $b_q(\alpha,\beta)=b_q(\beta,\beta)=1$, and 
$b_q(\beta,\xi)=b_q(\xi,\beta)=0$ otherwise. 
The condition $q=(S_q,\;\pi_q,\;b_q,\;u_q)$ is stronger than
$p\restriction\eta$, has $S_q\subseteq\eta$ and forces $\beta\in B_{\alpha}$,
$\beta\notin B_{\alpha_1},\dots,\beta\notin B_{\alpha_k}$. By Lemma 3.5
there is a condition $r$ that is stronger than both $p$ and $q$.
\enddemo

\medpagebreak

This concludes the proof that the forcing from Definition 3.3 adjoins a 
structure described in Theorem 3.1. That such a structure exists in $V$
is a consequence of the general theorem (Theorem 1.9) in \cite7. Our
forcing is $\omega_1$-uniform in the sense of Definition 1.1 in \cite7
and the dense sets needed to produce the $B_\alpha$ and the $A_n$ in
Theorem 3.1 conform to Definition 1.4 in \cite7
and hence the method of \cite7 applies.

\enddocument

\end